\documentclass[10pt,oneside,reqno,a4paper]{amsart}
\usepackage{setspace}
\usepackage{graphicx}
\usepackage{comment}
\usepackage{bbm}
\usepackage[dvipsnames]{xcolor}
\usepackage{multicol}
\usepackage[english]{babel}
\usepackage{amsmath,amssymb,graphicx,xcolor}
\usepackage[paperheight=297mm,paperwidth=220mm,margin=1in]{geometry}
\usepackage{times}
\usepackage{hyperref}
\hypersetup{
 colorlinks,
 citecolor=Maroon,
 linkcolor=Maroon,
 urlcolor=Maroon}

\newcommand{\assign}{:=}
\newcommand{\mathd}{\mathrm{d}}
\newcommand{\of}{:}
\newcommand{\tmabbr}[1]{#1}
\newcommand{\tmcolor}[2]{{\color{#1}{#2}}}
\newcommand{\tmem}[1]{{\em #1\/}}
\newcommand{\tmmathbf}[1]{\ensuremath{\boldsymbol{#1}}}
\newcommand{\tmname}[1]{\textsc{#1}}
\newcommand{\tmop}[1]{\ensuremath{\operatorname{#1}}}

\newcommand{\tmstrong}[1]{\textbf{#1}}

\newcommand{\tmtextit}[1]{{\itshape{#1}}}

\newcommand{\tmtextsc}[1]{{\scshape{#1}}}

\newenvironment{itemizedot}{\begin{itemize} }{\end{itemize}}

\newcommand{\zz}{\zeta}

\newcommand{\rev}[1]{#1}

\newcommand{\sh}{\mathcal{M}}
\newcommand{\drag}{\nu}
\newcommand{\eqass}{\tmcolor{black}{\assign}}
\newcommand{\sprthree}{\tmstrong{s{\tmname{p}}r\ensuremath{_3}}}
\newcommand{\II}{\mathcal{I}}
\newcommand{\XX}{\mathcal{X}}

\newcommand{\WW}{\mathcal{W}}

\newcommand{\YY}{\mathcal{Y}}

\newcommand{\LL}{\mathcal{L}}

\newcommand{\RR}{\mathbbmss{R}}
\newcommand{\NN}{\mathbbmss{N}}
\newcommand{\T}{\mathsf{T}}

\newcommand{\eqs}{=}

\begin{document}

\title{Parking 3-sphere swimmer\\\small{II.\,The {\rev{long arm}} asymptotic
  regime}}

\author{Fran{\c c}ois Alouges}
\address{\small Fran{\c c}ois Alouges,\, CMAP, Centre de Math{\'e}matiques Appliqu{\'e}es 
{\'E}cole Polytechnique\\
Route de Saclay, 91128 Palaiseau Cedex, France}

\author{Giovanni Di Fratta}
\address{\small Giovanni Di Fratta,\, Institute for Analysis and Scientific Computing, TU
Wien \\ Wiedner Hauptstra{\ss}e 8-10\\ 1040 Wien, Austria.}

\maketitle

      {\tmname{Abstract.}} The paper {\rev{carries on}} our previous
      investigations on the complementary version of Purcell's rotator
      ({\sprthree}): a low-Reynolds-number swimmer composed of three balls of
      equal radii. In the asymptotic regime of very long arms, the Stokes
      induced governing dynamics is derived, and then experimented in the
      context of energy minimizing self-propulsion characterized in the first
      part of the paper.
      
      \smallskip
     

\begin{multicols}{2}
\section{Introduction}

{\noindent}In his seminal paper~{\cite{purcell1977life}}, Purcell explains how
at small Reynolds numbers any organism trying to swim using the reciprocal
stroke of a scallop, which moves by opening and closing its valves, is
condemned to {\rev{go back to its original position at the end of one cycle}}.
This observation leads to the question of finding the simplest mechanisms
capable of self-propulsion at these scales; by this, we mean the ability to
moving by performing a cyclic shape change, a {\tmem{stroke}}, in the absence
of external forces. Several proposals have been put forward and analyzed (see,
e.g.,
{\cite{avron2004optimal,becker2003self,lefebvre2009stokesian,najafi2004simple,purcell1977life,taylor1951analysis}}
and the review paper {\cite{lauga2009hydrodynamics}}).

\begin{figure*}
  \raisebox{-0.482761964970761\height}{\includegraphics{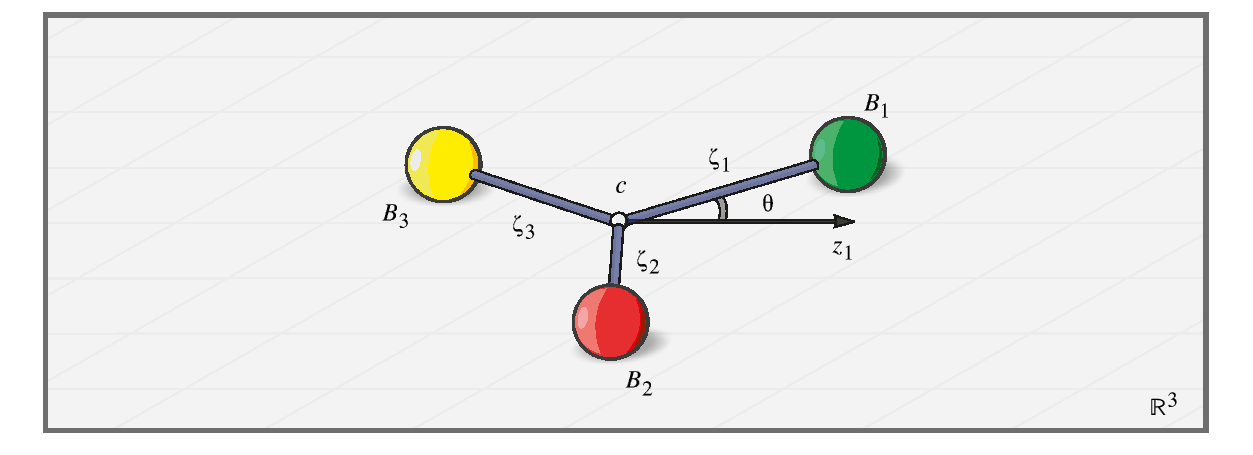}}
  \caption{\label{fig:1}The swimmer {\sprthree} \ is composed of three
  non-intersecting balls $(B_i)_{i \in \NN_3}$ of $\RR^3$ of equal radii. The
  three balls are linked together by {\tmem{thin arms}} that are able to
  elongate, independently of each other, along three coplanar axes that meet
  at the center $c \in \RR^3$ and make fixed angles of $2 \pi / 3$ from one to
  another.}
\end{figure*}

In this paper, we focus on a very {\rev{specific}} microswimmer: the
complementary version of Purcell's three-sphere rotator ({\sprthree})
introduced in {\cite{lefebvre2009stokesian}} and fully described in
Section~\ref{sec:descriptionspr3}. This swimmer consists of three
non-intersecting balls $(B_i)_{i \in \NN_3}$ of $\RR^3$ centered at $b_i \in
\RR^3$ and of equal radii $a > 0$ (for $n \in \NN$ we set $\NN_n \assign \{ 1,
\ldots, n \}$). The three balls can move along three coplanar axes that
mutually meet at a point $c \in \RR^3$, the {\tmem{center}}, with fixed angles
of $2 \pi / 3$ one to another; this reflects a situation where the balls are
linked together by very {\tmem{thin}} telescopic arms that can elongate (see
Figure~\ref{fig:1}). The swimmer can freely rotate around $c$ in the
horizontal plane containing the {\rev{arms}}, although owing to the symmetries
of the system, it is forced to stay in this plane.

Full controllability of {\sprthree}, as well as for a {\rev{broader}} class of
model swimmers, i.e., the ability of the swimmer to reach any point in the
plane with any orientation, has been proved in {\cite{alouges2013optimally}},
while analytical investigations on the optimal control problem have been the
object of {\cite{alouges2016spr31}}. Compared to
{\cite{alouges2013optimally,alouges2016spr31}}, we propose here a quantitative
analysis and we want to stress that, by contrast to the earlier works on the
topic ({\tmabbr{cf.}}
{\cite{agostinelli2018peristaltic,alouges2019energy,alouges2015can,avron2004optimal,becker2003self,dreyfus2005purcell,giraldi2015optimal,najafi2004simple,purcell1977life,tam2007optimal}}),
here the presence of three control variables and three position variables
makes the analysis more involved and rich.

The main aim of this second part is to put the optimality results proved in
{\cite{alouges2016spr31}} into a concrete setting, specifically: the Stokes
induced governing control system for {\sprthree} in the asymptotic regime of
very long arms. {\rev{First, we derive}} closed-form expressions for the
{\rev{dynamics}} and use asymptotic analysis to simplify the results. Then, we
focus on the analysis of energy-minimizing strokes, and we {\rev{identify}}
the{\rev{ optimal }}parameters of the control system in terms of the initial
length of the arms and the radius of the three balls. Finally, {\rev{we
present }}numerical simulations that show the qualitative features of the
optimal swimming style.

\section{Kinematics and Dynamics of {\sprthree}}
\label{sec:descriptionspr3}
\begin{figure*}
  \raisebox{-0.481916769662075\height}{\includegraphics{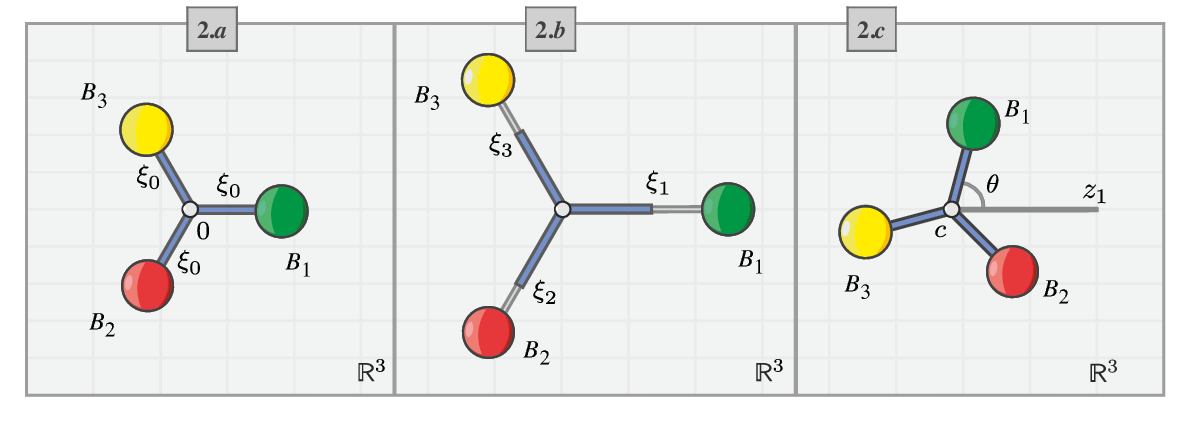}}
  \caption{\label{fig:j2}The swimmer {\sprthree} is fully described by the set
  of shape variables $\zz \eqass \left( \zz_1, \zz_2, \zz_3 \right) \in \sh$
  and by the set of position variables $p = (c, \theta) \in \RR^2 \times \RR$.
  [From left to right] In {\tmstrong{2.\,{\tmem{a}}}}, the reference
  configuration. The three spheres are located at the vertices of an
  equilateral triangle having the origin as barycenter ($c = 0$). In
  {\tmstrong{2.{\tmem{\,b}}}}, the set of shape variables $\zz \eqass (\xi_0 +
  \xi_1, \xi_0 + \xi_2, \xi_0 + \xi_3) \in \sh$ represents a possible shape
  state of the swimmer characterized by three different lengths of the arms.
  In {\tmstrong{2.\,{\tmem{c}}}}, a possible position state $(c, \theta)$,
  with $\theta \neq 0$, is sketched.}
\end{figure*}

{\noindent}As the three balls are not allowed to
rotate around their axes, the {\tmem{shape}} of the swimmer can be
parametrized by the lengths $\zz_1, \zz_2, \zz_3$ of its three arms, measured
from $c$ to the center of each of the balls. Therefore, the possible
geometrical configurations of the swimmer can be described by introducing two
sets of variables:
\begin{itemizedot}
  \item The vector of {\tmem{shape variables}} $\zz \eqass \left( \zz_1,
  \zz_2, \zz_3 \right) \in \sh \eqass \left( 2 a / \sqrt{3}, \infty \right)^3
 \subseteq \RR_+^3$
   from which relative distances $(b_{i_{} j})_{i, j \in
  \NN_3}$ between the balls are obtained ({\tmabbr{{\tmabbr{cf}}.}}
  {\eqref{eq:bieq}} and {\eqref{eq:conb}}). The lower bound on $\sh$ is
  imposed to exclude any overlap of the spheres.
  
  \item Position and orientation of $\sprthree$ in the plane are specified by
  the coordinates of the center $c \in \RR^2 \times \{ 0 \}$, and by the angle
  $\theta$ that one arm, {\tmabbr{e.g.}}, the arm connected to $B_1$, makes
  with the fixed direction $z_1$. We refer to $p = (c, \theta) \in \RR^2
  \times \RR$ as the vector of {\tmem{position variables}}.
\end{itemizedot}
Precisely, without loss of generality, we assume that in its initial
configuration, the three arms of the swimmer sit in the plane $\RR^2 \times \{
0 \}$. In order to compute the position of the three balls, we take the
vertices of the equilateral triangle defined as the convex hull of the unit
vectors $z_1, z_2, z_3 \in \RR^3$, with $z_1 \assign (1, 0, 0)^{\T}$, $z_2
\assign R^{\T} (2 \pi / 3) z_1$, $z_3 \assign R (2 \pi / 3) z_1$, and $R
(\phi)$ the planar rotation through an angle $\phi \in \RR$ around the vector
$\hat{e}_3 = (0, 0, 1)$:
\begin{equation}
  R (\phi) \eqass \left(\begin{array}{ccc}
    \cos \phi & - \sin \phi & 0\\
    \sin \phi & \cos \phi & 0\\
    0 & 0 & 1
  \end{array}\right) . \label{eq:Rotmat}
\end{equation}
Then, the center $b_i$ of the $i$-th ball of the swimmer is at position
({\tmabbr{cf.}} Figure \ref{fig:j2})
\begin{equation}
  b_i \assign c + \zz_i R (\theta) z_i \in \RR^2, \label{eq:bieq}
\end{equation}
where, here and in the sequel, we identify $\RR^2$ with $\RR^2 \times \{ 0 \}$
and, similarly, we identify the action of $R$ in {\eqref{eq:Rotmat}} with its
two-dimensional analog. Since the balls cannot intersect, the matrix $b \eqass
(b_1, b_2, b_3) \in \RR^{2 \times 3}$ is constrained to take values into the
set
\begin{equation}
  \mathcal{B} \assign \left\{ b \in \RR^{2 \times 3} \of \, \min_{i < j \,} |
  b_{i_{} j} | > 2 a_{_{_{_{_{}}}}} \right\}, \qquad b_{i_{} j} \assign b_i -
  b_j . \label{eq:conb}
\end{equation}
The time evolution of the swimmer can be traced through the state variables
$\left( \zz, p \right) \in \sh \times \RR^3$. For $i \in \NN_3$, the
instantaneous velocity of the $i$-th sphere is obtained by differentiating
relation {\eqref{eq:bieq}} with respect to time
\begin{equation}
  u_i \left( \zz, p \right) = \dot{c} + \dot{\zz}_i R (\theta) z_i +
  \dot{\theta} \zz_i R (\theta) z_i^{\bot}, \label{eq:uiclassicform}
\end{equation}
with $z_i^{\bot} \eqass R (\pi / 2) z_i$.

The viscous resistance of the arms is deemed negligible and,
therefore, we assume that the fluid fills up the whole space outside the
balls, that is, the exterior domain $\Omega \eqass \RR^3 \backslash \cup_{i =
1}^3 \bar{B}_i$. The geometry of $\Omega$ is uniquely determined by the common
radius $a$ of the three spheres, and by the matrix $b = (b_i)_{i \in \NN_3}$
having as columns the centers of the balls. At low Reynolds numbers, the
dynamics of the swimmer is governed by the Stokes equations
\begin{equation}
  \left\{ \begin{array}{rll}
    - \mu \Delta u + \nabla p \eqs & 0 \hspace{1.2em} & \text{in } \Omega,\\
    \tmop{div} u \eqs & 0 & \text{in } \Omega,
  \end{array} \right. \label{eq:Stokeseqs}
\end{equation}
where $u$ and $p$ are, respectively, the velocity field and the pressure of
the fluid, and $\mu$ is {\rev{its}} viscosity. As the structure of the swimmer
is {\rev{deformable but made of rigid balls}}, the governing equations are
subject to {\tmem{no-slip}} boundary conditions on the balls. Because of the
linearity of Stokes equations, the vector $\tmmathbf{u} \assign (u_1, u_2,
u_3)$ {\rev{collecting the three velocities}} in {\eqref{eq:uiclassicform}}
can be expressed in the algebraic form ({\tmabbr{cf.}}
{\cite{happel2012low,najafi2004simple}})
\begin{equation}
  \tmmathbf{u}=\mathcal{H}\tmmathbf{f}, \label{eq:HOseen}
\end{equation}
where $\mathcal{H}$ is the Oseen tensor (which depends on the viscosity) and
$\tmmathbf{f} \assign \rev{(f_1, f_2, f_3) \in \RR^6}$ is the vector
collecting the forces acting on the balls. Symmetry arguments show that in the
long arm asymptotic regime the hydrodynamic relation {\eqref{eq:HOseen}} takes
the form:
\begin{equation}
  u_i \eqs \frac{1}{\nu} f_i \; + \sum_{j \neq i \in \NN_3} \mathcal{S}
  (b_{i_{} j}) f_j \label{eq:ithvf}
\end{equation}
{\rev{where the {\tmem{stokeslet}}}}
\begin{equation}
  \mathcal{S} (x) \eqass \frac{1}{8 \pi \mu} \left( \frac{I}{| x |} + \frac{x
  \otimes x}{| x |^3} \right)
\end{equation}
represents a fundamental solution of the Stokes
system~{\cite{hancock1953self}}, and $\drag \assign 6 \pi \mu a \in \RR^+$ is
the drag coefficient linking, at small Reynolds numbers, the force to the
velocity of a spherical object of radius $a \in \RR^+$ immersed in a fluid of
viscosity $\mu$.

It will be convenient to rewrite {\eqref{eq:ithvf}} in the form
\begin{equation}
  \tmmathbf{u}= \left( \frac{1}{\drag} \II + \LL \right) \tmmathbf{f},
  \label{eq:Tbalgebraic}
\end{equation}
where $\II \assign \tmop{diag} (I, I, I)$ is the $6 \times 6$ identity matrix,
and $\LL$ the {\tmem{mutual interaction matrix}} defined by
\begin{equation}
  \LL \assign \left(\begin{array}{ccc}
    0 & \mathcal{S} (b_{12}) & \mathcal{S} (b_{13})\\
    \mathcal{S} (b_{12}) & 0 & \mathcal{S} (b_{23})\\
    \mathcal{S} (b_{13}) & \mathcal{S} (b_{23}) & 0
  \end{array}\right) . \label{eq:reluf1}
\end{equation}

\section{Dynamics of {\sprthree} in the limit of very long arms}

{\noindent}Due to the negligible inertia, the total viscous force and torque
exerted by the surrounding fluid on the swimmer must vanish. {\rev{In other
words, the dynamics is subject to the balance equations}}
\begin{equation}
  \sum_{i \in \NN_3} f_i = 0 \quad \text{and} \quad \sum_{i \in \NN_3} b_i
  \times f_i \eqs 0. \label{eq:totalfandtorque}
\end{equation}
{\rev{Here, the cross product}} stands for the determinant form on $\RR^2$ and
the $b_i$'s are given by {\eqref{eq:bieq}}. Clearly, for every $i \in \NN_3$,
there exist vectors $b_{\bot, i} (\zeta_i, \theta) \in \RR^2$, such that
$b_{\bot, i} (\zeta_i, \theta) \cdot f_i \eqs b_i \times f_i$, and, therefore,
the balance equations {\eqref{eq:totalfandtorque}} can be expressed in the
concise form
\begin{equation}
  \WW (\zeta, \theta) \tmmathbf{f}= 0, \label{eq:balanceeqsW}
\end{equation}
where the matrix $\mathcal{W}$ is defined by
\begin{equation}
  \WW (\zeta, \theta) \assign \left(\begin{array}{ccc} 
    I_{2 \times 2} & I_{2 \times 2} & I_{2 \times 2}\\[2pt]
    b_{\bot, 1}^{\T} (\zeta_1, \theta) & b_{\bot, 2}^{\T} (\zeta_2, \theta) &
    b_{\bot, 3}^{\T} (\zeta_3, \theta)
  \end{array}\right) \label{eq:torquematrix} .
\end{equation}
\tmcolor{red}{}{\rev{We assume that the three arms of the swimmer have the
same initial length $\xi_0 \in \RR^+$ with $\xi_0 \gg a$, and we set $\zz_i
\assign \xi_0 + \xi_i$ with $| \xi_i | \ll \xi_0$. We want to show that in the
limit of very long arms, and at the leading order, the swimming problem for
{\sprthree} reduces to a control problem of the form}}
\begin{equation}
  \dot{p} \; \eqs F (\xi, \theta) \dot{\xi}, \label{eq:controlsystemsp3}
\end{equation}
with $\xi \assign (\xi_1, \xi_2, \xi_3)$, {\rev{whose structural symmetries
have been fully investigated in the first part of the paper
({\tmabbr{cf.}}~{\cite{alouges2016spr31}}).}}

First, since the vector of the velocities depends linearly both on $\dot{\xi}$
and $\dot{p}$, we can recast relations {\eqref{eq:uiclassicform}} in the form
\begin{equation}
  \tmmathbf{u} \eqs \XX (\theta) \dot{\xi} + \YY (\xi, \theta) \dot{p} 
  \label{eq:eqforumatrix}
\end{equation}
where $\XX, \YY$ are the shape matrices given by
\begin{equation}
  \begin{array}{rll}
    \XX (\theta) & \assign & \left(\begin{array}{c|c|c}
      R (\theta) z_1 & 0_{2 \times 1} & 0_{2 \times 1}\\
      0_{2 \times 1} & R (\theta) z_2 & 0_{2 \times 1}\\
      0_{2 \times 1} & 0_{2 \times 1} & R (\theta) z_3
    \end{array}\right)_{6 \times 3},\\[8pt]\\
    \YY (\xi, \theta) & \assign & \left(\begin{array}{c|c}
      I_{2 \times 2 \;} & (\xi_0 + \xi_1) R (\theta) z_1^{\bot}\\
      I_{2 \times 2 \;} & (\xi_0 + \xi_2) R (\theta) z_2^{\bot}\\
      I_{2 \times 2 \;} & (\xi_0 + \xi_3) R (\theta) z_3^{\bot}
    \end{array}\right)_{6 \times 3} .
  \end{array} \label{eq:matricesPandX}
\end{equation}
In the limit of large arms, the mutual interaction matrix becomes a
perturbation of the diagonal part $(1 / \nu) \II$ and equation
{\eqref{eq:Tbalgebraic}} can be inverted to give (at the leading order)
\begin{equation}
  \tmmathbf{f}= \left( \drag \II - \drag^2 \LL \right) \tmmathbf{u} \eqs
  \left( \drag \II - \drag^2 \LL \right) \left( \XX (\theta) \dot{\xi} + \YY
  (\xi, \theta) \dot{p} \right) \label{eq:seventeen}
\end{equation}
by use of {\eqref{eq:eqforumatrix}}. Multiplying both members by
$\mathcal{W}$, and after simplifying by $\nu$, we infer that ({\tmabbr{cf.}}
{\eqref{eq:balanceeqsW}})
\begin{equation}
  \WW (\xi, \theta) \left( \II - \drag \LL \right) \left( \XX (\theta)
  \dot{\xi} + \YY (\xi, \theta) \dot{p} \right) \eqs 0,
\end{equation}
with the convenient and not dangerous abuse of notation $\mathcal{W} (\xi,
\theta) \assign \mathcal{W} (\zeta, \theta)$.
This is of the desired form {\eqref{eq:controlsystemsp3}} with
\begin{equation}
  F = - \left( \WW \left( \II - \drag \LL \right) \YY \right)^{- 1}  \WW 
  \left( \II - \drag \LL \right) \XX \label{eq:nineteen}
\end{equation}
where, to shorten notation, we left understood the parameters $\xi_0, \xi$ and
$\theta$. Moreover, because of the invariance of Stokes equations under the
group of rotations, according to {\cite{alouges2016spr31}} (Prop.~1), we can
factorize the control system $F$ in the form $F (\xi, \theta) \eqs R (\theta)
F (\xi)$ with $F (\xi) \assign F (\xi, 0)$, and therefore
\begin{equation}
  \dot{p} \; \eqs \; R (\theta) F (\xi) \dot{\xi} .
  \label{eq:control1spr3factorized}
\end{equation}
Also ({\tmabbr{cf.}}~{\cite{alouges2016spr31}} (Prop.~4)), in the limit of
small strokes, {\rev{i.e., in the regime $| \xi | / \xi_0 < a / \xi_0 \ll 1$
(see also Section \ref{sec:concremarks}),}} {\rev{we can expand $F$ to leading
order in $\xi$. This gives}}
\begin{equation}
  F (\xi) \dot{\xi} \eqs F_0 \dot{\xi} + \sum_{k \in \NN_3} (A_k \dot{\xi}
  \cdot \xi) e_k . \label{eq:control1spr3firstorder}
\end{equation}
Here, a straightforward computation shows that $F_0 \assign F (0)$ is given by
\begin{equation}
  F_0 = \varphi (a, \xi_0) \left( \begin{array}{ccc}
    - 2 & 1 & 1\\
    0 & \sqrt{3} & - \sqrt{3}\\
    0 & 0 & 0
  \end{array} \right),
\end{equation}
with
\begin{equation}
  \varphi (a, \xi_0) \assign \frac{1}{6} - \frac{1}{16 \sqrt{3}} (a / \xi_0)
  +\mathcal{O} (a / \xi_0)^2 .
\end{equation}
Instead, the first order correctors $(A_k)_{k \in \NN_3}$
({\tmabbr{cf.}}~{\cite{alouges2016spr31}} (Corollary~1)) have a special
structure which can be fully characterized in terms of four real parameters.
Precisely, there exist $\alpha \eqass \alpha (a, \xi_0), \beta \eqass \beta
(a, \xi_0), \gamma \eqass \gamma (a, \xi_0)$, and $\lambda \eqass \lambda (a,
\xi_0)$, depending only on the radius $a$ of the balls and on the initial
common length $\xi_0$ of the arms of {\sprthree}, such that
\begin{equation}
  \begin{array}{lll}
    A_1 & = & \left( \begin{array}{ccc}
      - \lambda & \alpha + \frac{1}{3} \beta & \alpha + \frac{1}{3} \beta\\
      - \alpha + \frac{1}{3} \beta & \frac{\lambda}{2} & - \frac{2}{3} \beta\\
      - \alpha + \frac{1}{3} \beta & - \frac{2}{3} \beta & \frac{\lambda}{2}
    \end{array} \right),\\ \\
    A_2 & = & \sqrt{3}  \left( \begin{array}{ccc}
      0 & \frac{\alpha - \beta}{3} & \frac{\beta - \alpha}{3}\\
      \frac{- \beta - \alpha}{3} & \frac{\lambda}{2} & - \frac{2 \alpha}{3}\\
      \frac{\alpha + \beta}{3} & \frac{2 \alpha}{3} & - \frac{\lambda}{2}
    \end{array} \right),
  \end{array} \label{eq:H12}
\end{equation}
and
\begin{equation}
  A_3 \eqs \left( \begin{array}{ccc}
    0 & - \gamma & \gamma\\
    \gamma & 0 & - \gamma\\
    - \gamma & \gamma & 0
  \end{array} \right) . \label{eq:H3}
\end{equation}
With the aid of a symbolic computation software and expanding $F$ in terms of
$\xi$ in {\eqref{eq:nineteen}} we can identify the entries and get
\begin{eqnarray}
  \alpha (a, \xi_0) & \assign & \frac{1}{\xi_0} \left( \frac{1}{32 \sqrt{3} }
  (a / \xi_0) +\mathcal{O} (a / \xi_0)^2 \right), \\
  \beta (a, \xi_0) & \assign & \frac{1}{\xi_0} \left( \frac{1}{16 \sqrt{3} }
  (a / \xi_0) +\mathcal{O} (a / \xi_0)^2 \right), \\
  \lambda (a, \xi_0) & \assign & \frac{1}{\xi_0} \left( \frac{5}{48 \sqrt{3}}
  (a / \xi_0) +\mathcal{O} (a / \xi_0)^2  \right), 
\end{eqnarray}
and
\begin{eqnarray}
  \gamma (a, \xi_0) & \assign & \frac{1}{\xi_0^2} \left( \frac{1}{6 \sqrt{3} }
  +\mathcal{O} (a / \xi_0)^2 \right) . 
\end{eqnarray}
We remark that only the skew-symmetric parts $(M_k)_{k \in
\NN_3}$ of the matrices $(A_k)_{k \in \NN_3}$ contribute to a net displacement
of the swimmer after one stroke ({\tmabbr{cf.}} {\cite{alouges2016spr31}}).
For any $\xi \in \RR^3$ they can be expressed by the actions
\begin{equation}
  M_1 \xi \eqs \alpha \xi \times \tau_1, \quad M_2 \xi \eqs \alpha \xi \times
  \tau_2, \quad M_3 \xi \eqs \gamma \xi \times \tau_3, \label{eq:actionMs}
\end{equation}
with
\begin{equation}
  \tau_1 \eqass (0, - 1, 1), \; \tau_2 \eqass \frac{1}{\sqrt{3}} (- 2, 1, 1),
  \; \tau_3 \eqass (1, 1, 1) \label{eq:taui}
\end{equation}
forming an orthogonal basis of $\RR^3$.

\section{Optimal swimming}
{\noindent}\label{sec:4}Following the notion of swimming efficiency proposed
by Lighthill in {\cite{lighthill1952squirming}} (cf.~also
{\cite{giraldi2015optimal,tam2007optimal}}), we adopt the following notion of
{\tmem{kinematic}} optimality: energy minimizing strokes are those minimizing
the kinetic energy dissipated during one stroke in order to reach a prescribed
net displacement $\delta p \in \RR^3$. In mathematical terms, the total energy
dissipation due to a smooth stroke $\zeta : I \rightarrow \sh$, can be
evaluated by considering the instantaneous power dissipated at time $t \in I$,
defined by $\mathcal{P} (\tmmathbf{u}) =\tmmathbf{f} \cdot \tmmathbf{u}$. We
note that $\dot{p}$ is linear in $\dot{\xi}$ because of
{\eqref{eq:controlsystemsp3}}, and so are $\tmmathbf{f}$ and $\tmmathbf{u}$
due to {\eqref{eq:eqforumatrix}} and {\eqref{eq:seventeen}}. Thus $\mathcal{P}
(\tmmathbf{u})$ turns out to be a quadratic form in $\dot{\xi}$ that we write
in the following form
\begin{equation}
  \mathcal{P} (\tmmathbf{u}) \eqs G (\xi) \dot{\xi} \cdot \dot{\xi}
\end{equation}
for a suitable matrix-valued function $G$ that, by the rotational invariance
of the problem, does not depend on $\theta$.

At the leading order in the limit of small strokes ({\tmabbr{cf.}}
{\cite{alouges2016spr31}} ({\textsection}$\,$5)) the instantaneous power
dissipated at time $t \in I$ reads as $\mathcal{P} (\tmmathbf{u} (t)) \eqs G_0
\dot{\xi} (t) \cdot \dot{\xi} (t)$, with $G_0 \assign G (0)$, and the total
energy dissipation associated with a stroke $\zeta : I \rightarrow \sh$ is
given by (recalling that $\zz_i \assign \xi_0 + \xi_i$)
\begin{equation}
  \mathcal{G} (\zeta) \eqass \int_I G_0  \dot{\xi} (t) \cdot \dot{\xi} (t)
  \mathd t . \label{eq:functionalenergydissipated}
\end{equation}
It can be readily checked that, as derived in {\cite{alouges2016spr31}}
({\textsection}$\,$5), the matrix $G_0$ is symmetric, positive-definite, and
has the following special structure
\begin{equation}
  G_0 = \left(\begin{array}{ccc}
    \kappa & h & h\\
    h & \kappa & h\\
    h & h & \kappa
  \end{array}\right),
\end{equation}
with the two parameters $h, \kappa$ depending only on the ratio $a / \xi_0$
between the radius of the balls of {\sprthree}, and on the common initial
length of its arms. Again, a symbolic computation shows that
\begin{equation}
  \begin{array}{lll}
    \kappa & \assign & \frac{2}{3} + \frac{1}{\sqrt{3}} (a / \xi_0)
    +\mathcal{O} (a / \xi_0)^2,\\ \\
    h & \assign & \frac{1}{6} + \frac{7}{16 \sqrt{3}} (a / \xi_0) +\mathcal{O}
    (a / \xi_0)^2 .
  \end{array} \label{eq:GMatrix}
\end{equation}
It is convenient to denote by $g_1 \eqass (\kappa - h)$ and $g_2 \eqass
(\kappa + 2 h)$ the eigenvalues of $G_0$. Note that $g_1$ is of multiplicity
two. Their expanded expressions read as
\begin{eqnarray}
  g_1 & \eqs &  \frac{1}{2} + \frac{3 \sqrt{3} }{16} (a / \xi_0) +\mathcal{O}
  (a / \xi_0)^2,  \label{eq:g1}\\
  g_2 & \eqs & 1 + \frac{5 \sqrt{3}}{8} (a / \xi_0) .  \label{eq:g3}
\end{eqnarray}
In {\cite{alouges2016spr31}} (Theorem~5.1)  we proved that the stroke $\zeta :
I \rightarrow \sh$ that produces a prescribed change of position and
orientation $\delta p \in \RR^3$ of the swimmer at the minimal cost
$\mathcal{G} (\zeta)$ is an ellipse of $\RR^3$. This optimal stroke is given
by
\begin{equation}
  \xi (t) \eqass (\cos t) u + (\sin t) v, \label{eq:ellipseF}
\end{equation}
where the vectors $u, v \in \RR^3$ can be fully computed from $\delta p$, the
coefficients $\alpha, \gamma$ of the skew-symmetric matrices $(M_k)_{k \in
\NN_3}$, and the eigenvalues of $G_0$.

Namely, as shown in {\cite{alouges2016spr31}} (Theorem~5.1), any minimizer
is, in $\xi$, an ellipse of $\RR^3$ centered at the origin, and the minimum
value of $\mathcal{G}$ is equal to $| \omega |$ where
\begin{equation}
  \omega \assign \tmop{diag} \left( \frac{\sqrt{g_1 g_2}}{\sqrt{2} \alpha},
  \frac{\sqrt{g_1 g_2}}{\sqrt{2} \alpha}, \frac{g_1 }{\sqrt{3} \gamma} \right)
  \delta p.
\end{equation}
More precisely, considering two orthogonal vectors $\rev{\varsigma_1},
\varsigma_2 \in \RR^3$ in the plane orthogonal to $\omega$ and such that $|
\rev{\varsigma_1} |^2 = | \rev{\varsigma_2} |^2 = | \omega |$, we can compute
the vectors $u$ and $v$ in {\eqref{eq:ellipseF}} via the relations
\begin{equation}
  u \eqass \frac{U \Lambda^{- 1 / 2}}{\sqrt{2 \pi}}  \rev{\varsigma_1}
  \hspace{0.17em}, \quad v \eqass \frac{U \Lambda^{- 1 / 2}}{\sqrt{2 \pi}}
  \varsigma_2, \label{eq:tocomputea1b1thm}
\end{equation}
with $U = (\tau_i / | \tau_i |)_{i \in \NN_3}$ ({\tmabbr{cf}}.
{\eqref{eq:taui}}) and $\Lambda \assign \tmop{diag} (g_1, g_1, g_2)$.

Summarizing, at the leading order in the range of small strokes and very long
arms, the governing dynamics of {\sprthree} for energy minimizing strokes is
given by ({\tmabbr{cf.}} {\eqref{eq:control1spr3firstorder}})
\begin{equation}
  \theta (t) \eqs \sigma t \quad \text{with} \quad \sigma \eqass \gamma (u
  \times v) \cdot \tau_3 \in \RR, \label{eq:ggd1}
\end{equation}
\begin{equation}
  \dot{c} (t) \eqs R (\sigma t)  F_0 \dot{\xi} (t) + R (\sigma t)
  \sum_{j \in \NN_2} (A_j  \dot{\xi} (t) \cdot \xi (t)) e_j, \label{eq:gd2}
\end{equation}
with $(A_j)_{j \in \NN_2}$ given by {\eqref{eq:H12}}, and $u, v \in \RR^3$
given by {\eqref{eq:tocomputea1b1thm}}. In particular, the angular velocity of
the swimmer is constant in time and is zero when the prescribed net
displacement $\delta p$ is purely translational ($\delta p_3 = 0$).

It is easily seen that energy minimizing net displacements along the $x$-axis
direction are achieved via elliptic strokes contained in the plane orthogonal
to the vector $\tau_1$. Similarly pure along-$y$ (resp. along-$\theta$) net
displacements are achieved via elliptic strokes contained in the plane
orthogonal to $\tau_2$ (resp. to $\tau_3$).

\begin{figure*}[t]
  \raisebox{-0.483312005341271\height}{\includegraphics{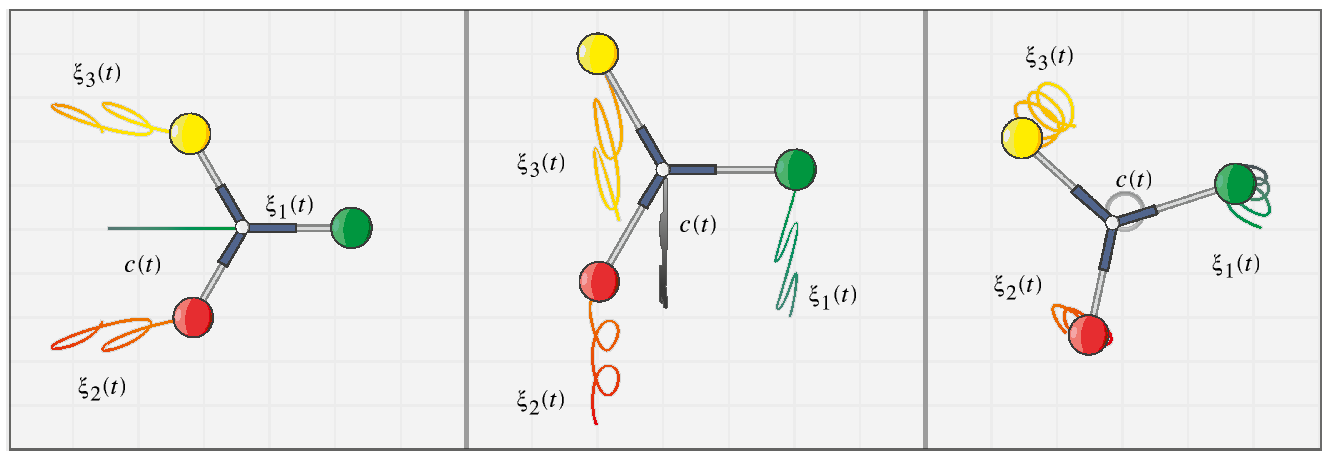}}
  \caption{\label{fig3}The evolution of $\xi$ and $c$ during optimal strokes.
  {\tmstrong{(left)}} Shape and position changes described by {\sprthree} to
  achieve a pure $x$ displacement. {\tmstrong{(center)}} Shape and position
  changes to achieve a pure $y$ displacement. {\tmstrong{(right)}} Shape and
  position changes to achieve a pure $\theta$ displacement.}
\end{figure*}

The results of numerical simulations of {\eqref{eq:ggd1}}-{\eqref{eq:gd2}}
when the control $\zeta$ is the optimal swimming strategy for a prescribed net
displacement $\delta p$ along the $x, y$ and $\theta$ directions are shown in
Figure~\ref{fig3}. Although we put some effort in drawing pictures that give a
good feeling of how the swimmer performs, we are aware that the dynamics can
be better appreciated by watching a video rather than looking at static
frames; in that regard, in the supplementary electronic material, it is
possible to find a video demonstrating the motion traced by {\sprthree} during
optimal swimming.
\section{Concluding Remarks}\label{sec:concremarks}
{\noindent}Note that, for $j = 1, 2$, we have that $\lim_{\xi_0 \rightarrow \infty} A_j (a,
\xi_0) =0$ and $ \lim_{a \rightarrow 0} A_j (a, \xi_0) \eqs 0$. However, since $\gamma
(0, \xi_0) = 1 / \left( 6 \sqrt{3} \xi_0^2 \right)$ and $A_3 = M_3$,
\begin{equation}
  \lim_{\xi_0 \rightarrow \infty} M_3 (a, \xi_0) \eqs 0, \quad \quad \lim_{a
  \rightarrow 0} M_3 (a, \xi_0) \neq 0 . \label{eq:asymmetry}
\end{equation}
In other words, the asymptotic limit of very small balls differs from one of
very long arms. This is understood by the presence of two fundamental
geometric scales: the common radius $a$ of the three balls, and the initial
length $\xi_0$ of its arms. In this respect, the two following asymptotic
regimes are different:
\begin{equation}
  \frac{a}{\xi_0} \ll \frac{| \xi |}{\xi_0} \ll 1, \quad \quad \frac{| \xi
  |}{\xi_0} \ll \frac{a}{\xi_0} \ll 1, \label{eq:asymmetryasympt}
\end{equation}
where we have denoted by $| \xi |$ the ``average'' stroke intensity.
\begin{itemizedot}
  \item In the limit $a / \xi_0 \ll | \xi | / \xi_0 \ll 1$ the swimmer offers
  great resistance to a net displacement in the $(x, y)$ coordinates, but it
  is strikingly still able to produce net angular displacements in the
  $\theta$ variable.
  
  \item The second condition in {\eqref{eq:asymmetryasympt}} represents the
  limit of very long arms and is more interesting for the applications as it
  allows for both translations and rotations.
\end{itemizedot}

\section{Acknowledgments}

{\noindent}This work was partially supported by the Labex LMH (grant
ANR-11-LABX-0056-LMH) in the {\tmem{Programme des Investissements d'Avenir}}.
Also, the second author acknowledges support from the Austrian Science Fund
(FWF) through the special research program {\tmem{Taming complexity in partial
differential systems}} (Grant SFB F65).

\end{multicols}
\end{document}